\documentclass[twoside,12pt,leqno]{article}
-\headheight\advance\topmargin by -\headsep

\usepackage{amsmath}
\usepackage{amsthm}
\usepackage{amssymb}
\usepackage{amscd}
\usepackage{graphicx}
\usepackage{epsfig}

\def\Z{\bf \mbox{Z\hspace{-.40em}Z}}
\def\R{{\bf \mbox{I\hspace{-.20em}R}}}

\def\lcf{\lbrack\! \lbrack}
\def\rcf{\rbrack\! \rbrack}
\numberwithin{equation}{section} \allowdisplaybreaks

\newtheorem{theorem}{Theorem}[section]

\newtheorem{proposition}{Proposition}[section]

\theoremstyle{definition}
\newtheorem{definition}{Definition}[section]

\begin{document}


\setcounter{page}{1} \thispagestyle{empty}



\markboth{F. Petalidou and J.M. Nunes da Costa}{Dirac structures
for generalized Courant and Courant algebroids}

\label{firstpage} $ $
\bigskip

\bigskip

\centerline{{\Large Dirac structures for generalized Courant and
Courant algebroids}}

\bigskip
\bigskip
\centerline{{\large by  Fani Petalidou and Joana M. Nunes da
Costa}}

\vspace*{.7cm}

\begin{abstract}
We establish some fundamental relations between Dirac subbundles
$L$ for the generalized Courant algebroid $(A\oplus A^{\ast},
\phi+W)$ over a differentiable manifold $M$ and the associated
Dirac subbubndles $\tilde{L}$ for the corresponding Courant
algebroid $\tilde{A} \oplus \tilde{A}^{\ast}$ over $M\times \R$.
\end{abstract}

\pagestyle{myheadings}
\section{Introduction}
In \cite{c}, T. Courant introduces the notion of a \emph{Dirac
structure} in order to present a unified framework for the study
of symplectic and Poisson structures and foliations. Alan
Weinstein and his collaborators develop the theory of these
structures and study several problems of Poisson geometry via
Dirac structures theory \cite{lwx1}, \cite{lwx2}. The notion was
exploited by A. Wade (\cite{wd}) and recently by the second author
and J. Clemente-Gallardo (\cite{jj}) in order to interpreter
Jacobi manifolds (\cite{lch}, \cite{dlm}) by means of Dirac
structures. In \cite{jj}, J.M. Nunes da Costa and J.
Clemente-Gallardo approach this problem by introducing the notions
of a \emph{generalized Courant algebroid} and of a \emph{Dirac
structure for a generalized Courant algebroid} and by proving that
the double $(A\oplus A^{\ast},\phi+W)$ of a generalized Lie
bialgebroid $((A,\phi),(A^{\ast},W))$ over a differentiable
manifold $M$, notion very close to the Jacobi manifolds
(\cite{im1}), is a generalized Courant algebroid.

In the present work, being well known that there is an one-to-one
correspondence between generalized Lie bialgebroids structures
$((A,\phi),(A^{\ast},W))$ over $M$ and Lie bialgebroids structures
$(\tilde{A},\tilde{A}^{\ast})$, $\tilde{A}=A\times \R$,
$\tilde{A}^{\ast}=A^{\ast}\times \R$, over $\tilde{M}=M\times \R$,
we establish some basic relations between the Dirac subbundles $L$
for $(A\oplus A^{\ast}, \phi+W)$ and the associated Dirac
subbundles $\tilde{L} = \{X+e^t\alpha \, /\, X+ \alpha \in L\}$
for $\tilde{A} \oplus \tilde{A}^{\ast}$. We prove : 1) $L$ is a
reducible Dirac structure for $(A\oplus A^{\ast},\phi+W)$ if and
only if $\tilde{L}$ is a reducible Dirac structure for $\tilde{A}
\oplus \tilde{A}^{\ast}$. 2) If $\mathcal{F}$ and
$\tilde{\mathcal{F}}$ are the characteristic foliations of $M$ and
$\tilde{M}$ defined by $L$ and $\tilde{L}$, respectively, then
$\tilde{L}$ induces an homogeneous Poisson structure on
$\tilde{M}/\tilde{\mathcal{F}}=M/\mathcal{F}\times \R$ which is
the Poissonization of the induced Jacobi structure on
$M/\mathcal{F}$ by $L$.

\vspace{2mm} \noindent \textbf{Notation :} In this paper, $M$ is a
$C^{\infty}$-differential manifold of finite dimension. We denote
by $C^{\infty}(M)$ the space of all real
$C^{\infty}$-differentiable functions on $M$ and by $\delta$ the
usual de Rham differential operator.

\section{Generalized Lie bialgebroids}
Let $(A, [\,,\,], a)$ be a Lie algebroid over $M$ (\cite{mck}),
$A^{\ast}$ its dual vector bundle over $M$, $\bigwedge A^{\ast} =
\oplus_{k \in \Z}\bigwedge^k A^{\ast}$ the graded exterior algebra
of $A^{\ast}$ whose differential sections are called {\it A-forms}
on $M$, $d : \Gamma(\bigwedge A^{\ast}) \to \Gamma(\bigwedge
A^{\ast})$ the exterior derivative of degree 1 and $\phi \in
\Gamma(A^{\ast})$ an 1-cocycle in the Lie algebroid cohomology
complex with trivial coefficients (\cite{mck}, \cite{im1}), i.e.,
for all $X,Y\in \Gamma(A)$, $\langle \phi, [X,Y]\rangle =
a(X)(\langle \phi, Y\rangle) - a(Y)(\langle \phi, X\rangle)$. We
modify the usual representation $a$ of the Lie algebra
$(\Gamma(A), [\,,\,])$ on the space $C^{\infty}(M)$ by defining
$a^{\phi} : \Gamma(A)\times C^{\infty}(M)\to C^{\infty}(M)$,
$a^{\phi}(X,f) = a(X)f + \langle \phi,X\rangle f$. The resulting
cohomology operator $d^{\phi} : \Gamma(\bigwedge A^{\ast}) \to
\Gamma(\bigwedge A^{\ast})$ of the new cohomology complex is
called the $\phi$-{\it differential} of $A$ and $d^{\phi}\eta =
d\eta + \phi\wedge \eta$, for all $\eta \in \Gamma(\bigwedge ^k
A^{\ast})$. $d^{\phi}$ allows us to define the $\phi$-{\it Lie
derivative by} $X\in \Gamma(A)$, $\mathcal{L}_X^{\phi} :
\Gamma(\bigwedge^k A^{\ast}) \to \Gamma(\bigwedge^k A^{\ast})$, as
$\mathcal{L}_X^{\phi} = d^{\phi}\circ i_X + i_X \circ d^{\phi}$,
where $i_X$ is the contraction by $X$. Using $\phi$ we can also
modify the Schouten bracket $[\,,\,]$ on $\Gamma(\bigwedge A)$ to
the $\phi$-{\it Schouten bracket} $[\,,\,]^{\phi}$ on
$\Gamma(\bigwedge A)$ by setting, for all $P\in \Gamma(\bigwedge^p
A)$ and $Q \in \Gamma(\bigwedge^q A)$, $[P,Q]^{\phi} = [P,Q] +
(p-1)P\wedge (i_{\phi}Q) + (-1)^p(q-1)(i_{\phi}P)\wedge Q$, where
$i_{\phi}Q$ can be interpreted as the usual contraction of a
multivector field with an 1-form. For details, see \cite{mck},
\cite{im1} and \cite{gm1}.

The notion of \emph{generalized Lie bialgebroid} has been
introduced by D. Iglesias and J.C. Marrero in \cite{im1} and
independently by J. Grabowski and G. Marmo in \cite{gm1} under the
name of {\it Jacobi bialgebroid}, in such a way that a Jacobi
manifold (\cite{lch}) has a generalized Lie bialgebroid
canonically associated and conversely. We recall that a
\emph{Jacobi manifold} is a smooth manifold $M$ equipped with a
bivector field $\Lambda$ and a vector field $E$ such that
$[\Lambda,\Lambda]=-2E \wedge \Lambda$ and $[E,\Lambda] = 0$,
where $[\,,\,]$ denotes the Schouten bracket.

We consider a Lie algebroid $(A, [\,,\,], a)$ over $M$ and an
1-cocycle $\phi \in \Gamma(A^{\ast})$ and we assume that the dual
vector bundle $A^{\ast}\to M$ admits a Lie algebroid structure
$([\,,\,]_{\ast}, a_{\ast})$ and that $W \in \Gamma(A)$ is an
1-cocycle in the Lie algebroid cohomology complex with trivial
coefficients of $(A^{\ast},[\,,\,]_{\ast}, a_{\ast})$. Then, we
say that :
\begin{definition}\label{D1}
The pair $((A,\phi),(A^{\ast},W))$ is a \emph{generalized Lie
bialgebroid} over $M$ if, for all $X,Y\in \Gamma(A)$ and $P\in
\Gamma(\bigwedge^pA)$, the following conditions hold :
$$
d_{\ast}^W[X,Y] = [d_{\ast}^WX,Y]^{\phi} + [X,d_{\ast}^WY]^{\phi}
\hspace{3mm} \mathrm{and} \hspace{3mm} \mathcal{L}_{\ast \phi}^WP
+ \mathcal{L}_W^{\phi}P = 0 \,;
$$
$d_{\ast}^W$ and $\mathcal{L}_{\ast}^W$ are, respectively, the
$W$-differential and the $W$-Lie derivative of $A^{\ast}$.
\end {definition}

Obviously, if $\phi = 0$ and $W=0$, we recover the notion of {\it
Lie bialgebroid} introduced by K. Mackenzie and P. Xu in \cite{mx}
and its equivalent definition given by Yv. Kosmann-Schwarzbach in
\cite{ks1}.

\vspace{2mm}

Given a Lie algebroid $(A,[\,,\,],a)$ over $M$, we can construct a
Lie algebroid structure on $\tilde{A} \to \tilde{M}$,
$\tilde{A}=A\times \R$ and $\tilde{M}=M\times \R$. We identify
$\Gamma(\tilde{A})$ with the set of the time-dependent sections of
$A \to M$, i.e. for any $\tilde{X}\in \Gamma(\tilde{A})$ and
$(x,t)\in M\times \R$, $t$ being the canonical coordinate on $\R$,
$\tilde{X}(x,t) = \tilde{X}_t(x)$, where $\tilde{X}_t \in
\Gamma(A)$, and we take : i) the Lie bracket
$[\,,\,]^{\,\tilde{}}$ on $\Gamma(\tilde{A})$ defined, for any
$\tilde{X},\tilde{Y}\in \Gamma(\tilde{A})$ and $(x,t) \in
\tilde{M}$, by $[\tilde{X},\tilde{Y}]^{\,\tilde{}}(x,t) =
[\tilde{X}_t,\tilde{Y}_t](x)$, ii) the bundle map $\tilde{a} :
\tilde{A}\to T\tilde{M}$, $\tilde{a}(\tilde{X})(x,t) =
a(\tilde{X}_t)(x)$. Then $(\tilde{A}, [\,,\,]^{\,\tilde{}},
\tilde{a})\to \tilde{M}$ is a Lie algebroid. Also, taking an
1-cocycle $\phi$ of $A$, we deform
$([\,,\,]^{\,\tilde{}},\tilde{a})$ in two different ways and we
obtain two new Lie algebroid structures on $\tilde{A}$,
\cite{im1}. Precisely, for any $\tilde{X},\tilde{Y}\in
\Gamma(\tilde{A})$:
\begin{equation}\label{11}
[\tilde{X},\tilde{Y}]^{\,\tilde{}\phi} =
[\tilde{X},\tilde{Y}]^{\,\tilde{}} + i_{\phi}\tilde{X}_t
\partial\tilde{Y}/\partial t - i_{\phi}\tilde{Y}\partial\tilde{X}/\partial t,
\hspace{3mm}\tilde{a}^{\phi}(\tilde{X}) = \tilde{a}(\tilde{X})
+i_{\phi}\tilde{X}\partial/\partial t;
\end{equation}
\begin{eqnarray}\label{12}
[\tilde{X},\tilde{Y}]^{\,\hat{}\phi} & =&  e^{-t}(
[\tilde{X},\tilde{Y}]^{\,\tilde{}} + \langle
\phi,\tilde{X}_t\rangle (\partial\tilde{Y}/\partial t-\tilde{Y}) -
\langle \phi,\tilde{Y}_t\rangle
(\partial \tilde{X}/\partial t-\tilde{X})),\\
\hat{a}^{\phi}(\tilde{X}) & = &  e^{-t}( \tilde{a}(\tilde{X}) +
\langle \phi,\tilde{X}_t\rangle \partial/\partial t). \nonumber
\end{eqnarray}

\begin{theorem}[\cite{im1}]\label{genbial-bial}
Let $(A,[\,,\,],a)$ be a Lie algebroid over $M$ and $\phi \in
\Gamma(A^{\ast})$ an 1-cocycle. Suppose that $A^{\ast}$ has a Lie
algebroid structure $([\,,\,]_{\ast},a_{\ast})$ and that $W\in
\Gamma(A)$ is an 1-cocycle for this structure. Consider on
$\tilde{A} = A\times \R$ and $\tilde{A}^{\ast} = A^{\ast}\times
\R$ the Lie algebroid structures $([\,,\,]^{\,\tilde{}\phi},
\tilde{a}^{\phi})$ and $([\,,\,]^{\,\hat{}\,W}_{\ast},
\hat{a}_{\ast}^W)$, respectively. Then
$(\tilde{A},\tilde{A}^{\ast})$ is a Lie bialgebroid over
$\tilde{M}=M\times \R$ if and only if $((A,\phi),(A^{\ast},W))$ is
a generalized Lie bialgebroid over $M$. The induced Poisson
structure on $\tilde{M}$ is the Poissonization of the induced
Jacobi structure on $M$.
\end{theorem}

Moreover, the image $\mathrm{Im}a$ of the anchor map $a$ of
$(A,[\,,\,],a)\to M$ is an integrable distribution on $M$
(\cite{df}) which defines a singular foliation $\mathcal{F}_A$ of
$M$, called the {\it Lie algebroid foliation of M associated with
A} (\cite{im2}). The relation between the leaves of the Lie
algebroid foliation $\mathcal{F}_{\tilde{A}}$ of $M\times \R$
associated with
$(\tilde{A},[\,,\,]^{\,\tilde{}\phi},\tilde{a}^{\phi})$ (given by
(\ref{11})) and the leaves of the Lie algebroid foliation
$\mathcal{F}_A$ of $M$ associated with $A$ was studied in
\cite{im2} by D. Iglesias and J.C. Marrero. They have proved :

\begin{theorem}[\cite{im2}]\label{foliat}
Under the above considerations, suppose that $(x_0,t_0)\in M\times
\R$ and that $\tilde{F}$ and $F$ are the leaves of the Lie
algebroid foliations $\mathcal{F}_{\tilde{A}}$ and $\mathcal{F}_A$
passing through $(x_0,t_0)\in M\times \R$ and $x_0\in M$,
respectively, and denote by $A_{x_0}$ the fiber of $A$ over $x_0$.
Then : (1) If $\ker(a\vert_{A_{x_0}})\nsubseteq
\langle\phi(x_0)\rangle^{\circ}$, $\tilde{F} = F\times \R$. (2) If
$\ker(a\vert_{A_{x_0}})\subseteq \langle\phi(x_0)\rangle^{\circ}$
and $\pi_1 : M\times \R \to M$ is the canonical projection onto
the first factor, $\pi_1(\tilde{F}) = F$ and
$\pi_1\vert_{\tilde{F}} : \tilde{F} \to F$ is a covering map.
\end{theorem}

\section{Generalized Courant algebroids}
The notion of {\it generalized Courant algebroid} has been
introduced by the second author and J. Clemente-Gallardo in
\cite{jj} and independently, under the name of {\it Courant-Jacobi
algebroid}, by J. Grabowski and G. Marmo in \cite{gm2}.

\begin{definition}[\cite{jj}]\label{Def-GenCour}
Let $E\to M$ to be a vector bundle over a differentiable manifold
$M$ equipped with : (i) a nondegenerate symmetric bilinear form
$(\,,\,)$ on the bundle, (ii) a skew-symmetric bilinear bracket
$[\,,\,]$ on $\Gamma(E)$, (iii) a bundle map $\rho : E \to TM$ and
(iv) an $E$-1-form $\theta$ such that, for any $e_1,e_2 \in
\Gamma(E)$, $\langle\theta, [e_1,e_2]\rangle = \rho(e_1)\langle
\theta,e_2\rangle - \rho(e_2)\langle\theta,e_1\rangle$. We
consider : (a) the bundle map $\rho^{\theta} : E \to TM \times \R$
defined, for any $e \in E$, by $\rho^{\theta}(e) = \rho(e) +
\langle\theta,e\rangle$, (b) the applications $\mathcal{D},
\mathcal{D}^{\theta} : C^{\infty}(M) \to \Gamma(E)$ defined, for
any $f\in C^{\infty}(M)$, respectively, by $\mathcal{D}f =
\frac{1}{2}\beta^{-1}\rho^{*}\delta f$ \footnote{$\beta$ is the
isomorphism from $E$ onto $E^*$ given by the nondegenerate
bilinear form $(\,,\,)$.} and $\mathcal{D}^{\theta}f =
\mathcal{D}f + \frac{1}{2}f \beta^{-1}(\theta)$ and (c) for any
$e_1,e_2,e_3 \in \Gamma(E)$, the function $T(e_1,e_2,e_3)=
\frac{1}{3}([e_1,e_2],e_3) + c.p.$ on the base $M$. Then, we say
that $E$ is a \emph{generalized Courant
algebroid} if the following relations are satisfied : \\
1. $[[e_1,e_2],e_3] + c.p. = \mathcal{D}^{\theta}T(e_1,e_2,e_3),
\hspace{5mm} \forall \, e_1,e_2,e_3 \in \Gamma(E)$; \\
2. $\rho^{\theta}([e_1,e_2]) = [\rho^{\theta}(e_1),
\rho^{\theta}(e_2)], \footnote{The bracket on the right-hand side
is the Lie bracket defined on $\Gamma(TM\times \R)$ by
$[(X,f),(Y,g)]=([X,Y], X\cdot g - Y\cdot f)$.} \hspace{5mm}
\forall
\, e_1,e_2 \in \Gamma(E)$; \\
3. $[e_1,fe_2] = f[e_1, e_2] + (\rho(e_1)f)e_2 -
(e_1,e_2)\mathcal{D}f, \hspace{3mm} \forall \, e_1,e_2 \in
\Gamma(E), \;\forall f\in C^{\infty}(M)$; \\
4. $\rho^{\theta}\circ \mathcal{D}^{\theta} = 0$, i.e., for any
$f,g \in C^{\infty}(M)$,
$(\mathcal{D}^{\theta}f,\mathcal{D}^{\theta}g) = 0$;\\
5. $\rho^{\theta}(e)(e_1,e_2) = ([e,e_1]
+\mathcal{D}^{\theta}(e,e_1),e_2) + (e_1,[e,e_2]
+\mathcal{D}^{\theta}(e,e_2)), \hspace{1mm} \forall\, e,e_1,e_2
\in \Gamma(E)$.
\end{definition}

\begin{definition}\label{D2}
A \emph{Dirac structure for a generalized Courant algebroid}
$(E,\theta)$ over $M$ is a subbundle $L\subset E$ that is maximal
isotropic under $(\,,\,)$ and integrable, i.e. $\Gamma(L)$ is
closed under $[\,,\,]$.
\end{definition}

A Dirac subbundle $L$ of $(E,\theta)$ is a Lie algebroid under the
restrictions of the bracket $[\,,\,]$ and of the anchor $\rho$ to
$\Gamma(L)$. If $\theta \in \Gamma(L^{\ast})$, then it is an
1-cocycle for the Lie algebroid cohomology with trivial
coefficients of $(L,[\,,\,]\vert_{L},\rho\vert_L)$.

\vspace{2mm}

The most important example of generalized Courant algebroid is the
double $(A\oplus A^{\ast}, \phi+W)$ of a generalized Lie
bialgebroid $((A,\phi),(A^{\ast},W))$ over $M$. On $A\oplus
A^{\ast}$ there exist two natural nondegenerate bilinear forms,
one symmetric and another skew-symmetric $(\,,\,)_{\pm}$ : for any
$X_1+\alpha_1, X_2+\alpha_2 \in A\oplus A^{\ast}$,
$(X_1+\alpha_1,X_2+\alpha_2)_{\pm} = 1/2(\langle
\alpha_1,X_2\rangle \pm \langle \alpha_2,X_1\rangle)$ and on
$\Gamma(A\oplus A^{\ast})\cong \Gamma(A)\oplus \Gamma(A^{\ast})$
we introduce the bracket $\lcf \,,\,\rcf$ : for all $X_1 +
\alpha_1, X_2 + \alpha_2 \in \Gamma(A\oplus A^{\ast})$,
\begin{eqnarray*}\label{14}
\lcf X_1 + \alpha_1,X_2 + \alpha_2 \rcf & \!\!\! = \!\!\! &
([X_1,X_2]^{\phi} + \mathcal{L}_{\ast \alpha_1}^WX_2 -
                         \mathcal{L}_{\ast \alpha_2}^WX_1 - d_{\ast}^W(e_1,e_2)_- ) + \\
                  & \!\!\! + \!\!\! & ([\alpha_1,\alpha_2]_{\ast}^W + \mathcal{L}_{X_1}^{\phi}\alpha_2 -
                         \mathcal{L}_{X_2}^{\phi}\alpha_1 +
                         d^{\phi}(e_1,e_2)_-).
\end{eqnarray*}
Also, we consider the bundle map $\rho : A\oplus A^{\ast} \to TM$
given by $\rho = a + a_{\ast}$, i.e., for any $X+\alpha \in E$,
$\rho(X+\alpha) = a(X) + a_{\ast}(\alpha)$. We have :

\begin{theorem}[\cite{jj}]
If $((A,\phi),(A^{\ast},W))$ is a generalized Lie bialgebroid over
$M$, then $A \oplus A^{\ast}$ endowed with $(\lcf \,,\,\rcf,
(\,,\,)_+,\rho)$ and $\theta =\phi + W\in \Gamma(E^{\ast})$ is a
generalized Courant algebroid over $M$. The operators
$\mathcal{D}$ and $\mathcal{D}^{\theta}$ are, respectively,
$\mathcal{D} = (d_{\ast}+d)\vert_{C^{\infty}(M)}$ and
$\mathcal{D}^{\theta} = (d_{\ast}^W +
d^{\phi})\vert_{C^{\infty}(M)}$.
\end{theorem}

\section{Dirac structures of $((A,\phi),(A^{\ast},W))$ and of $(\tilde{A},\tilde{A}^{\ast})$}
Let $((A,[\,,\,],a,\phi),(A^{\ast},[\,,\,]_{\ast},a_{\ast},W))$ be
a generalized Lie bialgebroid over $M$ and $(A\oplus A^{\ast},
\lcf \,,\,\rcf,(\,,\,)_+, a + a_{\ast},\phi+W)$ the associated
generalized Courant algebroid.

\begin{definition}\label{Def-red}
We say that a Dirac subbundle $L$ of $A\oplus A^{\ast}$ is
\emph{reducible} if the image $a(D)$ of its \emph{characteristic
subbundle} $D=L\cap A$ by $a$ defines a simple foliation
$\mathcal{F}$ of $M$. By the term "simple foliation" we mean that
$\mathcal{F}$ is a regular foliation such that the space
$M/\mathcal{F}$ is a nice manifold and the canonical projection $M
\to M/\mathcal{F}$ is a submersion.
\end{definition}

\begin{definition}\label{Def-admis}
Let $L$ be a Dirac subbundle of $A\oplus A^{\ast}$. A function $f
\in C^{\infty}(M)$ is called \emph{$L$-admissible} if there exists
$Y_f \in \Gamma(A)$ such that $Y_f + d^{\phi}f \in \Gamma(L)$. We
denote by $C_L^{\infty}(M,\R)$ the set of all $L$-admissible
functions of $C^{\infty}(M)$.
\end{definition}

Let $((\tilde{A},[\,,\,]^{\,\tilde{}\phi},\tilde{a}^{\phi}),
(\tilde{A}^{\ast},[\,,\,]_{\ast}^{\,\hat{}\,W},\hat{a}_{\ast}^W))$
be the Lie bialgebroid over $\tilde{M}$ defined by
$((A,[\,,\,],a,\phi),(A^{\ast},[\,,\,]_{\ast},a_{\ast},W))$ as in
Theorem \ref{genbial-bial}. Then,
$\tilde{A}\oplus\tilde{A}^{\ast}$ endowed with : (i) the two
nondegenerate bilinear forms $(\,,\,)_{\pm}$ on
$\tilde{A}\oplus\tilde{A}^{\ast}$ : for all
$\tilde{X}_1+\tilde{\alpha}_1, \tilde{X}_2+\tilde{\alpha}_2 \in
\tilde{A}\oplus\tilde{A}^{\ast}$,
$(\tilde{X}_1+\tilde{\alpha}_1,\tilde{X}_2+\tilde{\alpha}_2)_{\pm}
= 1/2(\langle
\tilde{\alpha}_1,\tilde{X}_2\rangle\pm\langle\tilde{\alpha}_2,\tilde{X}_1\rangle)$,
(ii) the bracket $\lcf\,,\,\rcf^{\,\tilde{}}$ on
$\Gamma(\tilde{A}\oplus\tilde{A}^{\ast})$ : for all
$\tilde{X}_1+\tilde{\alpha}_1, \tilde{X}_2+\tilde{\alpha}_2 \in
\Gamma(\tilde{A}\oplus\tilde{A}^{\ast})$,
\begin{eqnarray*}
\lcf
\tilde{X}_1+\tilde{\alpha}_1,\tilde{X}_2+\tilde{\alpha}_2\rcf^{\,\tilde{}}
& = & \big([\tilde{X}_1,\tilde{X}_2]^{\,\tilde{}\phi} +
                         \hat{\mathcal{L}}_{\tilde{\alpha}_1}^W\tilde{X}_2-
                         \hat{\mathcal{L}}_{\tilde{\alpha}_2}^W\tilde{X}_1 -
                         \hat{d}_{\ast}^W((\tilde{e}_1,\tilde{e}_2)_{-})\big)+  \\
                  &  & \big([\tilde{\alpha}_1,\tilde{\alpha}_2]_{\ast}^{\,\hat{}\,W} +
                        \tilde{\mathcal{L}}_{\tilde{X}_1}^{\phi}\tilde{\alpha}_2-
                        \tilde{\mathcal{L}}_{\tilde{X}_2}^{\phi}\tilde{\alpha}_1 +
                        \tilde{d}^{\phi}((\tilde{e}_1,\tilde{e}_2)_{-})\big),
\end{eqnarray*}
(for any $\tilde{f}\in C^{\infty}(\tilde{M})$,
$\tilde{d}^{\phi}\tilde{f} = \tilde{d}\tilde{f}
+\frac{\partial\tilde{f}}{\partial t}\phi$ and
$\hat{d}_{\ast}^W\tilde{f}=e^{-t}(\tilde{d}\tilde{f}+\frac{\partial\tilde{f}}{\partial
t}\phi)$, \cite{im1}), (iii) the bundle map $\tilde{\rho} :
\tilde{A}\oplus\tilde{A}^{\ast} \to T\tilde{M}$, $\tilde{\rho} =
\tilde{a}^{\phi}+\hat{a}_{\ast}^W$, is a Courant algebroid over
$\tilde{M}$ (\cite{lwx1}).

\vspace{2mm}\noindent Let $\mathbf{E} : \Gamma(A\oplus A^{\ast})
\to \Gamma(\tilde{A}\oplus\tilde{A}^{\ast})$ be the embedding of
$\Gamma(A\oplus A^{\ast})$ into
$\Gamma(\tilde{A}\oplus\tilde{A}^{\ast})$ defined, for any
$X+\alpha \in \Gamma(A\oplus A^{\ast})$, by
$$
\mathbf{E}(X+\alpha) = X + e^t\alpha,
$$
where $X$ and $\alpha$ are regarded as time-independent sections
of $\tilde{A}$ and $\tilde{A}^{\ast}$, respectively. If $L$ is a
subbundle of $A \oplus A^{\ast}$, we write $\tilde{L} =
\mathbf{E}(L)$ in order to denote the vector subbundle $\tilde{L}$
of $\tilde{A}\oplus \tilde{A}^{\ast}$ whose space of global cross
sections is the image by $\mathbf{E}$ of the space of global cross
sections of $L$, i.e. $\Gamma(\tilde{L}) = \mathbf{E}(\Gamma(L))$.
\begin{proposition}
Let $L$ be a vector subbundle of $A\oplus A^{\ast}$ and $\tilde{L}
= \mathbf{E}(L)$. Then, $L$ is a Dirac structure for the
generalized Courant algebroid $(A\oplus A^{\ast}, \phi+W)$ if and
only if $\tilde{L}$ is a Dirac structure for the Courant algebroid
$\tilde{A}\oplus \tilde{A}^{\ast}$.
\end{proposition}
\begin{proof}
It is easy to check that $\tilde{L}$ is a maximally isotropic
subbundle of $(\tilde{A}\oplus \tilde{A}^{\ast},(\,,\,)_+)$ if and
only if $L$ is a maximally isotropic subbundle of $(A\oplus
A^{\ast},(\,,\,)_+)$. Moreover, by a straightforward calculation
we get that
$$
\lcf \mathbf{E}(X_1+\alpha_1),\mathbf{E}(X_2+\alpha_2)
\rcf^{\,\tilde{}} = \mathbf{E}(\lcf
X_1+\alpha_1,X_2+\alpha_2\rcf), \hspace{3mm} \forall \,
X_1+\alpha_1, X_2+\alpha_2\in \Gamma(L),
$$
i.e. $\Gamma(\tilde{L})$ is closed under $\lcf
\,,\,\rcf^{\,\tilde{}}$ if and only if $\Gamma(L)$ is closed under
$\lcf \,,\,\rcf$.
\end{proof}

\begin{proposition}\label{admis}
Let $L$ be a Dirac structure of $(A\oplus A^{\ast}, \phi+W)$ and
$\tilde{L}=\mathbf{E}(L)$ the associated Dirac structure of
$\tilde{A}\oplus \tilde{A}^{\ast}$. Then $\tilde{f}\in
C^{\infty}(\tilde{M})$ is a $\tilde{L}$-admissible function if and
only if $\tilde{f}=e^tf$ and $f\in C_L^{\infty}(M)$.
\end{proposition}
\begin{proof}Let $\tilde{f}\in C_{\tilde{L}}^{\infty}(\tilde{M})$, i.e. there exists $Y\in
\Gamma(A)$ : $Y +\tilde{d}^{\phi}\tilde{f} \in \Gamma(\tilde{L})$.
But, $Y +\tilde{d}^{\phi}\tilde{f} \in \Gamma(\tilde{L})$ implies
that there exists $\xi\in \Gamma(A^{\ast})$ : $Y +\xi \in
\Gamma(L)$ and $Y +\tilde{d}^{\phi}\tilde{f} = \mathbf{E}(Y
+\xi)$, thus $\tilde{d}^{\phi}\tilde{f} = e^t\xi$. From
\emph{Theorem of normal forms for Lie algebroids} (\cite{df}) we
have that, if the rank of $a(D)$, $D=L\cap A$, at a point $q\in M$
is $k$, then we can construct on a neighborhood $U$ of $q$ in $M$
a system of local coordinates $(x_1,\ldots,x_k,\ldots,x_n)$ ($n =
\dim M$) and a basis of sections $(X_1,\ldots,X_k,\ldots X_r)$ of
$\Gamma(A)$ ($r$ is the dimension of the fibres of $A\to M$), with
$(X_1,\ldots,X_k)$ sections of $\Gamma(D)$, such that $a(X_i) =
\frac{\partial}{\partial x_i}$, for every $i = 1,\ldots,k$. Let
$(\alpha_1,\ldots,\alpha_k,\ldots,\alpha_r)$ be the basis of
$\Gamma(A^{\ast})$, dual of $(X_1,\ldots,X_k,\ldots X_r)$. Since
$\phi,\xi \in \Gamma(A^{\ast})$, there exist $\phi_i,\xi_i\in
C^{\infty}(U)$, $i=1,\ldots,r$, such that $\phi =
\sum_{i=1}^r\phi_i\alpha_i$ and $\xi = \sum_{i=1}^r\xi_i\alpha_i$.
So, for any $i = 1,\ldots,r$,
\begin{equation}\label{22}
\tilde{d}^{\phi}\tilde{f} = e^t\xi \Rightarrow \langle
\tilde{d}\tilde{f}+ (\partial\tilde{f}/\partial t)\phi,X_i\rangle
= \langle e^t\xi,X_i\rangle \Leftrightarrow \langle
\tilde{d}\tilde{f},X_i\rangle + (\partial\tilde{f}/\partial
t)\phi_i = e^t\xi_i.
\end{equation}
But, for $i=1,\ldots,k$, $\langle \tilde{d}\tilde{f},X_i\rangle =
\langle \delta \tilde{f}, \tilde{a}(X_i)\rangle = \langle \delta
\tilde{f}, a(X_i)\rangle = \langle \delta \tilde{f},
\frac{\partial}{\partial x_i}\rangle = \frac{\partial
\tilde{f}}{\partial x_i}$. Hence, the last equation of (\ref{22})
can be written, for any $i=1,\ldots,k$, as
\begin{equation}\label{23}
\partial \tilde{f}/\partial x_i + (\partial
\tilde{f}/\partial t)\phi_i = e^t\xi_i.
\end{equation}
By resolving the characteristic system $\frac{\delta
x_i}{1}=\frac{\delta t}{\phi_i}=\frac{\delta \tilde{f}}{e^t\xi_i}$
of (\ref{23}), we obtain that $\tilde{f}$ must be, at least
locally, of the form $\tilde{f}=e^tf$ with $f\in C^{\infty}(U)$.
Taking into account Definition \ref{Def-admis} and that
$\tilde{L}=\mathbf{E}(L)$, we get $\tilde{f} = e^tf\in
C_{\tilde{L}}^{\infty}(\tilde{M}) \Leftrightarrow f\in
C_L^{\infty}(M)$.
\end{proof}

\begin{proposition}\label{Rel-red}
Let $L$ be a Dirac subbundle for $(A\oplus A^{\ast}, \phi+W)$ and
$\tilde{L}=\mathbf{E}(L)$ the associated Dirac subbundle of
$\tilde{A}\oplus \tilde{A}^{\ast}$. Then, $L$ is reducible if and
only if $\tilde{L}$ is reducible.
\end{proposition}
\begin{proof}
Let $D=L\cap A$ and $\tilde{D}=\tilde{L}\cap \tilde{A}$ be the
characteristic subbundles of $L$ and $\tilde{L}$, respectively,
$\mathcal{F}$ and $\tilde{\mathcal{F}}$ the foliations of $M$ and
$\tilde{M}$, respectively, defined by $a(D)$ and
$\tilde{a}^{\phi}(\tilde{D})$, respectively. Obviously,
$\tilde{D}\cong D$ and $\tilde{a}^{\phi}(\tilde{D}) =
\{\tilde{a}^{\phi}(X)\, / \, X \in D\} = \{a(X) + \langle \phi,
X\rangle \partial /\partial t \, /\, X\in D\}$. Let $(x_0,t_0)$ be
a point of $\tilde{M} = M\times \R$, $\tilde{F}$ and $F$ the
leaves of $\tilde{\mathcal{F}}$ and $\mathcal{F}$ passing through
$(x_0,t_0)\in \tilde{M}$ and $x_0 \in M$, respectively, and
$D_{x_0}$ the fibre of $D$ over $x_0$. By Theorem \ref{foliat}, we
have : (i) if $ker(a\vert_{D_{x_0}})\nsubseteq\langle
\phi(x_0)\rangle^{\circ}$, then $\tilde{F}=F\times \R$, so $\dim
\tilde{F} = \dim F + 1$ and the vector field $\partial / \partial
t$ is tangent to $\tilde{F}$; (ii) if
$ker(a\vert_{D_{x_0}})\subseteq\langle \phi(x_0)\rangle^{\circ}$
and $\pi_1 : M\times \R \to M$ is the canonical projection, then
$\pi_1(\tilde{F})=F$ and $\pi_1\vert_{\tilde{F}} : \tilde{F}\to F$
is a covering map, thus $\dim \tilde{F} = \dim F$ and the vector
field $\partial /\partial t$ is not tangent to $\tilde{F}$. Since
every $\tilde{L}$-admissible function $\tilde{f}$ is of type
$\tilde{f}=e^tf$, $f\in C_L^{\infty}(M)$, (Proposition
\ref{admis}) and also it is constant along the leaves of
$\tilde{\mathcal{F}}$ (\cite{c},\cite{lwx2}), it is not possible
the leaves $\tilde{F}$ of $\tilde{\mathcal{F}}$ to be of type
$\tilde{F}=F\times \R$ (because, in this case, $\partial /
\partial t$ is tangent to $\tilde{F}$ and $\tilde{f} = e^tf$ is
not constant along $\partial / \partial t$). Thus, for any leaf
$\tilde{F}$ of $\tilde{\mathcal{F}}$ and for the corresponding
leaf $F$ of $\mathcal{F}$, we have : $\pi_1(\tilde{F})=F$ and
$\pi_1\vert_{\tilde{F}} : \tilde{F}\to F$ is a covering map.
Hence, we get : (1) Every leaf $\tilde{F}$ of
$\tilde{\mathcal{F}}$ is of the same dimension as the
corresponding leaf $F$ of $\mathcal{F}$, so $\mathcal{F}$ is a
regular foliation of $M$ if and only if $\tilde{\mathcal{F}}$ is a
regular foliation of $\tilde{M}$. (2) $\tilde{\mathcal{F}}\cong
\mathcal{F}$, so $\tilde{M} / \tilde{\mathcal{F}} \cong (M\times
\R) / \mathcal{F} \cong (M/\mathcal{F})\times \R$ ; thus,
$M/\mathcal{F}$ is a nice manifold if and only if $\tilde{M} /
\tilde{\mathcal{F}}$ is a nice manifold and the projection $M \to
M / \mathcal{F}$ is a submersion if and only if the projection
$M\times \R =\tilde{M} \to \tilde{M} / \tilde{\mathcal{F}} \cong
(M/\mathcal{F})\times \R$ is a submersion. Consequently, $L$ is a
reducible Dirac subbundle for $A\oplus A^{\ast}$ if and only if
$\tilde{L} = \mathbf{E}(L)$ is a reducible Dirac subbundle for
$\tilde{A}\oplus \tilde{A}^{\ast}$.
\end{proof}

Let $L$ be a Dirac structure of $(A\oplus A^{\ast},\phi+W)$ and
$\tilde{L}$ the associated Dirac structure of $\tilde{A}\oplus
\tilde{A}^{\ast}$. On $C_L^{\infty}(M)$ we define the bracket
$\{\,,\,\}_L$ by setting, for all $f,g \in C_L^{\infty}(M)$,
$\{f,g\}_L : = \rho^{\theta}(e_f)g$, where $e_f = Y_f + d^{\phi}f
\in \Gamma(L)$. Also, on $C_{\tilde{L}}^{\infty}(\tilde{M})$ we
define the bracket $\{\,,\,\}_{\tilde{L}}$ by setting, for all
$\tilde{f},\tilde{g}\in C_{\tilde{L}}^{\infty}(\tilde{M})$,
$\tilde{f}=e^tf$, $\tilde{g}=e^tg$ with $f,g\in C_L^{\infty}(M)$,
$\{\tilde{f},\tilde{g}\}_{\tilde{L}} : =
\tilde{\rho}(\tilde{e}_{\tilde{f}})\tilde{g}$, where
$\tilde{e}_{\tilde{f}} = Y_f + \tilde{d}^{\phi}\tilde{f} \in
\Gamma(\tilde{L})$. By a straightforward calculation we get :
\begin{equation}\label{eq-hom}
\{\tilde{f},\tilde{g}\}_{\tilde{L}} = \{e^tf,e^tg\}_{\tilde{L}} =
e^t\{f,g\}_L.
\end{equation}

\begin{theorem}[\cite{fj}]
1) If $1\in C_L^{\infty}(M)$\footnote{We have (\cite{fj}) : $1$ is
an $L$-admissible function if and only if, for any $Y\in
\Gamma(D)$, $\langle \phi,Y\rangle$ = 0.}, then
$(C_L^{\infty}(M),\{\,,\,\}_L)$ is a Jacobi algebra. 2) If $L$ is
a reducible Dirac subbundle of $(A\oplus A^{\ast},\phi+W)$ and
$1\in C_L^{\infty}(M)$, then $L$ induces a Jacobi structure on
$M/\mathcal{F}$ defined by the Jacobi bracket $\{\,,\,\}_L$.
\end{theorem}

\begin{theorem}
1) If $1\in C_L^{\infty}(M)$, then
$(C_{\tilde{L}}^{\infty}(\tilde{M}), \{\,,\,\}_{\tilde{L}})$ is an
homogeneous Poisson algebra with respect $\partial / \partial
t$\footnote{In the sense of Dazord-Lichnerowicz-Marle (\cite{dlm})
terminology.}. 2) If $L$ is a reducible Dirac subbundle of
$(A\oplus A^{\ast},\phi+W)$ and $1\in C_L^{\infty}(M)$, then
$\tilde{L}$ induces an homogeneous Poisson structure on
$\tilde{M}/\tilde{\mathcal{F}}$ defined by the homogeneous Poisson
bracket $\{\,,\,\}_{\tilde{L}}$. 3) $\tilde{M}/\tilde{\mathcal{F}}
= (M/\mathcal{F})\times \R$ and the induced homogeneous Poisson
structure on $\tilde{M}/\tilde{\mathcal{F}}$ by $\tilde{L}$ is the
Poissonization of the induced Jacobi structure on $M/\mathcal{F}$
by $L$.
\end{theorem}
\begin{proof}
1) It is checked by taking account (\ref{eq-hom}) and the fact
that, if $1\in C_L^{\infty}(M)$, then
$(C_L^{\infty}(M),\{\,,\,\}_L)$ is a Jacobi algebra. 2) By
applying the results of \cite{lwx2} to the reducible Dirac
subbundle $\tilde{L}$ and the homogeneous Poisson algebra
$(C_{\tilde{L}}^{\infty}(\tilde{M}), \{\,,\,\}_{\tilde{L}})$. 3)
We have $\tilde{\mathcal{F}} = \mathcal{F}\times \{0\}$
(\cite{fj}), thus $\tilde{M}/\tilde{\mathcal{F}} =
(M/\mathcal{F})\times \R$, and by (\ref{eq-hom}) we conclude the
announced result.
\end{proof}


\begin{tabular}{lcl}
Fani Petalidou &  &  Joana M. Nunes da Costa \\
{\it Faculty of Sciences and Technology}  & & {\it Department of Mathematics} \\
{\it University of Peloponnese} & & {\it University of Coimbra} \\
{\it 22100 Tripoli, Greece}  & & {\it Apartado 3008} \\
& & {\it 3001-454 Coimbra, Portugal} \\
 &  & \\
{\it e-mail : petalido@uop.gr} & & {\it e-mail :
jmcosta@mat.uc.pt}
\end{tabular}

\label{lastpage}
\end{document}